\documentclass{cmslatex}
\usepackage{amssymb,amsmath}
\usepackage{graphicx}%
\renewcommand{\theequation}{\arabic{section}.\arabic{equation}}

\newtheorem{thm}{Theorem}[section]
\newtheorem{lem}[thm]{Lemma}
\newtheorem{cor}[thm]{Corollary}

\newcommand{\R}{\mathbb R}

\newcommand{\teps}{\tau_{\eps}}
\newcommand{\reps}{\rho_{\eps}}
\newcommand{\meps}{{\cal M}^{\eps}}
\newcommand{\Teps}{T_{\reps} \meps}
\newcommand{\veps}{v_{\eps}}
\newcommand{\geps}{\gamma_{\eps}}

\newcommand{\ueps}{u_{\eps}}

\newcommand{\tveps}{{\tilde v}_{\eps}}
\newcommand{\jeps}{j_{\eps}}
\newcommand{\bigpar}{\par\smallskip}
\def\longrightharpoonup{
\relbar\joinrel\joinrel\relbar\joinrel\joinrel\relbar\joinrel\joinrel\rightharpoonup}
\newcommand{\xrightharpoonup}[1]{\stackrel{#1}{\longrightharpoonup}}
\newcommand{\Rset}{{\mathbb{R}}}
\newcommand{\pair}[2]{{\left({#1},\,{#2}\right)}}

\newcommand{\at}[1]{{\left({#1}\right)}}

\newcommand{\bat}[1]{{\big(#1\big)}}
\newcommand{\Bat}[1]{{\Big(#1\Big)}}
\newcommand{\norm}[1]{{\|{#1}\|}}
\newcommand{\abs}[1]{\left|{#1}\right|}
\newcommand{\babs}[1]{\big|{#1}\big|}

\newcommand{\eps}{{\varepsilon}}

\newcommand{\calM}{\mathcal{M}}
\newcommand{\calT}{\mathcal{T}}
\begin{document}
\title{%
Kramers' formula for chemical reactions in the context of Wasserstein gradient flows%
%
}%
\author{%
Michael Herrmann\thanks{%
Oxford Centre for Nonlinear PDE (OxPDE), (michael.herrmann@maths.ox.ac.uk).
}%
\and %
Barbara Niethammer\thanks {OxPDE, (niethammer@maths.ox.ac.uk).}%
}%
\pagestyle{myheadings}
\markboth{Kramers' formula in the context of Wasserstein gradient flows}{M. Herrmann and B.
Niethammer}
\maketitle%
\begin{center}
\today\\\smallskip
\end{center}

\begin{abstract}
We derive Kramers' formula as singular limit of the Fokker-Planck equation with double-well
potential. The convergence proof is based on the Rayleigh principle of the underlying
Wasserstein gradient structure  and complements a recent result by Peletier, Savar\'e and
Veneroni.
\end{abstract}
%
%
\section{Introduction}
%
%
In 1940 Kramers derived chemical reaction rates from certain limits in a Fokker-Planck 
equation that describes the probability density of a Brownian particle in an energy landscape
\cite{Kramers40}. The limit of high activation energy has been revisited in a recent paper by
Peletier et al \cite{PSV10}, where a spatially inhomogeneous extension of Kramers' formula is
rigorously derived for unimolecular reactions between two chemical states $A$ and $B$. Their
derivation relies on passing to the limit in the $L^2$-gradient flow structure of the
Fokker-Planck equation. It is well-known by now that the Fokker-Planck equation has also an
interpretation as a Wasserstein gradient flow \cite{JKO98} and the question was raised in
\cite{PSV10} whether Kramers' formula can also be derived and interpreted within this
Wasserstein gradient flow structure. This concept has also been investigated on a formal
level for more complicated reaction-diffusion systems in \cite{Mielke10}. A further
motivation for studying the Fokker-Planck equation within the Wasserstein framework comes
from applications that additionally prescribe the time evolution of a moment \cite{Dreyer10}.
\par
In this note we present a rigorous derivation which is based on passing to the limit within
the Wasserstein gradient flow structure. To keep things simple we restrict ourselves to the
spatially homogeneous case and consider the simplest case of a unimolecular reaction between
two chemical states $A$ and $B$, which are represented as two wells of an enthalpy function
$H:\Rset\to\Rset$. To avoid unimportant technicalities we assume that the enthalpy function
$H$ is a `typical' double-well potential. Specifically, we assume that $H$ is
a smooth, nonnegative and even function that satisfies
\begin{align*}
\text{$xH^\prime\at{x}<0\;$  for  $\;0<\abs{x}<1$},\qquad
\text{$xH^\prime\at{x}>0\;$  for  $\;\abs{x}>1$},
\end{align*}
with
\begin{align*}
H\at{\pm1}=H^\prime\at{\pm1}=H^\prime\at{0}=0,\quad H\at{0}=H\at{\pm2}=1, \quad
H^{\prime\prime}\at{0}<0<H^{\prime\prime}\at{\pm1}.
\end{align*}
The probability density of a molecule with chemical state $x$ is in the following denoted by
$\rho$. In Kramers' approach the molecule performs a Brownian motion in the energy landscape
described by $H$, so the evolution of $\rho$ is governed by the Kramers-Smoluchowski equation
\begin{align}\label{KSeq}
\partial_t \rho  = \partial_x \bat{\partial_x \rho + \eps^{-2}\rho H^\prime}\,,
\end{align}
where $\eps$ is the so-called 'viscosity' coefficient. In what follows we consider the high
activation energy limit $\eps^2 \ll 1$.
\par
The leading order dynamics of \eqref{KSeq} can be derived by formal asymptotics and governs
the evolution of
\begin{align*}
u^+(t) := 2\int_0^{\infty} \rho(t,x)\,dx\,,\qquad u^-(t):=2\int_{\infty}^0
\rho(t,x)\,dx\,,
\end{align*}
which satisfy $u^-\at{t}+u^+\at{t}=2$ due to $\int_\Rset{\rho}\pair{t}{x}\,dx=1$. Using WKB
methods, for instance, one finds
\begin{align}
\label{kramers}
\frac{d u^\pm}{dt}  = \frac{\tau_\eps k }{2}\big( u^\mp-u^\pm\big)
,\qquad%
\tau_\eps:=\eps^{-2} e^{-1/\eps^2}
,\qquad%
k:=  { \pi}^{-1} \sqrt{|H^{''}(0)| H^{''}(1)}.
\end{align}
We emphasize that the constant $k$ depends on the details of the function $H$ near its two
local minima and its local maximum, and that the time scale is exponentially slow in the
height barrier $H(0)/\eps^2 = 1/\eps^2$ between the two wells.
\par
Our goal is to derive \eqref{kramers} rigorously by passing to the limit $\eps \to 0$ in
\eqref{KSeq}. In order to derive a non-trivial limit we have to rescale time accordingly.
Thus we consider in the following the  probability distribution $\reps = \reps(t,x)$ that is
a solution of
\begin{align}
\label{eq1}
\teps\partial_t\reps =
\partial_x \bat{\partial_x\reps + \eps^{-2}\reps{H^\prime}}\,,
\qquad x \in \R, \;t>0.
\end{align}
An important role in the analysis will be played by the unique invariant measure
\begin{align}
\label{gammadef}
\geps(x) := Z_{\eps}^{-1} e^{-H(x)/\eps^2} \,, \qquad
Z_{\eps} := \int_{\R} e^{-H(y)/\eps^2} \,dy = \eps\frac{2\sqrt{2\pi}}{\sqrt{H^{''}(1)}}\big( 1+
o(1)\big),
\end{align}
which converges in the weak$\star$ topology of probability measures to
$\frac12(\delta_{-1}+\delta_{+1})$,
where $\delta_{\pm1}$ denotes the delta distributions in $\pm1$,
\par
As in \cite{PSV10} it is often convenient to switch to the density of $\reps$ with respect to
$\geps$, that is $\ueps=\reps/\geps$. Heuristically we expect that $\ueps$ is -- to leading
order in $\eps$ -- piecewise constant for $x<0$ and $x>0$, where the respective values
correspond to $u^+$ and $u^-$ as introduced above. In what follows we write
$u$ instead of $u^+$, so $u^-$ is given by $2-u$.
\par
For the derivation of the limit equation we assume that our
data are
well-prepared.
\bigpar%
\begin{thm}\label{T.1}
Let $\teps$ and $k$ be as in \eqref{kramers}, and for each $\eps$ let
$\reps:[0,\infty)\times\Rset\to[0,\infty)$ be a solution to
\eqref{KSeq} with initial datum $\reps^0$. Moreover, suppose that
the initial data $\reps^0$ are probability measures on $\Rset$ that
converge weakly$\star$ as $\eps\to0$ to some probability measure
$\rho^0$, and satisfy
\begin{align*}
\int_{\Rset} \geps|\ueps^0|^2\,dx
\,+\, \int_{\Rset} \frac{\geps|\partial_x \ueps^0|^2}{\teps} \,dx \leq C\,,
\qquad\ueps^0(x)= \reps^0(x)/\geps(x)\geq c>0\quad\forall\;x\in\R
\end{align*}
with constants $C$ and $c$ independent of $\eps$. Then, for all $t\geq 0$ we have that
\begin{align*}
\reps(t,\cdot)\xrightharpoonup{\eps\to0}\tfrac12\Big(u(t) \delta_{+1} + \big (2-u(t)\big)
\delta_{-1} \Big)
\end{align*}
weakly$\star$ in the space of probability measures, where the function
$u:[0,\infty)\to[0,\infty)$ satisfies
\begin{align}
\label{limiteq}
\dot u = -  k \at{u - 1},\qquad\qquad
u(0)=u_0\,,
\end{align}
where $u_0=2\int_0^\infty\rho^0\,dx$.%
\end{thm}
\bigpar
This result has already been derived in \cite{PSV10} in the more general setting with spatial
diffusion and under slightly weaker assumptions on the initial data. Our main contribution
here is therefore not the result as such, but the method of proof. We answer the question
posed in \cite{PSV10},  how the passage to the limit can be performed within a Wasserstein
gradient flow structure and we identify the corresponding structure for the limit.
\par
We present the formal gradient flow structures of \eqref{eq1} as well as \eqref{limiteq} in
section~\ref{sec:Structure}. In order to derive the limit equation we pass to the limit in
the Rayleigh principle that is associated to any gradient flow. This strategy is inspired by
the notion of $\Gamma$-convergence and has already been successfully employed in other
singular limits of gradient flows (e.g. in \cite{NO1}). In section~\ref{S.compactness} we
first obtain some basic a priori estimates as well as an approximation for $\ueps$ that will
be essential in the identification of the limit gradient flow structure. It is somewhat
unsatisfactory that we cannot derive suitable estimates solely from the energy estimates
associated to the Wasserstein gradient flow. Instead we use the estimates that correspond to
the $L^2$-gradient flow structure that is satisfied by $\ueps$. It is not obvious to us how
this can be avoided.
\par
Section~\ref{S.proofs} finally contains our main result, that is the novel proof of
Theorem~\ref{T.1}. We show that the limit of $\partial_t \reps$ satisfies the Rayleigh
principle that one obtains as a limit of the Rayleigh principle associated to the Wasserstein
gradient structure of \eqref{eq1}. As a consequence the limit is a solution of
\eqref{limiteq}.
%
%
%
\section{Gradient flows and Rayleigh principle}
\label{sec:Structure}
%
We briefly summarize the Wasserstein gradient structure of the Fokker-Planck equation as
well as the corresponding gradient flow structure of the limit problem. To point out the key
ideas we give a formal exposition and postpone some technical details to Section
\ref{S.compactness.1}.
\par
Given an energy functional $E$ on a manifold ${\cal M}$, whose tangent space $T_x{\cal M}$
is endowed with a metric tensor $g_x$, the $g$-gradient flow $x(t)$ of $E$ is defined such
that
\begin{align}\label{gradflow1}
g_{x(t)}\big( \dot x(t), v\big) + DE(x(t)) v =0
\end{align}
for all $v \in T_{x(t)}{\cal M}$ and for all $t>0$. Here $DE(x(t)) v$ denotes the directional
derivative of $E$ in direction $v$. Our convergence result relies on the Rayleigh principle,
that amounts to the observation that a curve $x(t)$ on $\calM$ solves \eqref{gradflow1} if
and only if for each $t$ the derivative $\dot x(t)$ minimizes
\begin{align*}
\tfrac 1 2 g_{x(t)} \big( v,v\big) + DE(x(t))v
\end{align*}
among all $v \in T_{x(t)}{\cal M}$. Slightly more general -- and more suitable for
generalizations to abstract manifolds in function spaces -- is the time integrated version:
For each $0<T<\infty$ the function $\dot{x}$ minimizes
\begin{align*}
\int_0^{T} \tfrac 1 2 g_{x(t)} \big( v(t),v(t)\big) + DE(x(t))v(t)\,dt
\end{align*}
among all functions $v$ with $v(t) \in T_{x(t)}{\cal M}$ for all $t \in [0,T]$. In what
follows we pass to the limit $\eps\to0$ in the time integrated Rayleigh principle because
only this one is compatible with the boundedness and compactness results derived below.
\par
To describe the {\bf{Wasserstein gradient structure}} of the $\eps$--problem \eqref{eq1} we
consider the formal manifold
\begin{align}
\notag
\meps:= \Big\{ \reps\colon \R \to [0,\infty) \,,\, \int_{\Rset}  \reps \,dx =1\Big\},\qquad
\Teps := \Big\{ \veps \colon \R \to \R \Big\}\,,
\end{align}
along with the metric tensor
\begin{align}
\label{eq12}
g^{\eps}_{\reps}\big( \veps,\veps\big) := \teps \int_{\Rset} \reps |\jeps|^2 \,dx,
\qquad \mbox{ where }  \quad
\veps=-\partial_x (\reps \jeps).
\end{align}
The energy of the Fokker-Planck equation is given by
\begin{align}
\notag
E^{\eps}(\reps) := \int_{\Rset} \reps \bat{\ln \reps +\eps^{-2}H}\,dx,
\end{align}
and has the directional derivative
\begin{align}
\label{eq14}
DE^{\eps}(\reps)\veps= \int_{\Rset} \veps
 \big( \ln \reps + \eps^{-2}H\big) \,dx = \int_{\Rset} \reps \jeps \partial_x\big
( \ln \reps + \eps^{-2}H\big)\,.
\end{align}
Consequently, the direction of steepest descent $\veps$ is characterized by the
requirement that
\begin{align}
\notag
g^{\eps}_{\reps}(\veps, \tveps) = - DE^{\eps}(\reps)\tveps \qquad \forall\;\tveps
 \in \Teps.
\end{align}
This means
\begin{align}
\notag
\jeps = -\teps^{-1}\partial_x\bat{\ln \reps + \eps^{-2}H}\,,
\end{align}
and we conclude that \eqref{eq1} is in fact the formal $g^{\eps}$-gradient flow of $E^{\eps}$
on $\meps$.
\par
The {\bf{gradient structure of the limit problem}} is very simple. The corresponding manifold
\begin{align}
\notag
{\cal M}:=\big\{ u \in (0,2)\big\}\,,\qquad
T_{u}{\cal M} := \big\{ v\in \R\big\}
\end{align}
is one-dimensional and equipped with the metric tensor
\begin{align}
\notag
g_{u}\big(v,v\big) := \frac{v^2\,\ln (u/(2-u))}{ 2 k (u-1)}\,.
\end{align}
Notice that the metric tensor is continuous in $u$ with $g_1\pair{v}{v}=v^2/k$, and
that $k >0$ is defined in
\eqref{kramers}. The limit energy is given by
\begin{align}
\notag
E(u) := \tfrac 1 2 \big (u \ln u + (2-u) \ln (2-u)\big)\,,\qquad
DE(u)v = \tfrac{1}{2}\,{v}\ln\bat{u/(2-u)}\,,
\end{align}
and we easily check that \eqref{limiteq} is the $g$-gradient flow of $E$ on ${\cal M}$.
%
%
\section{A priori estimates and implications}
\label{S.compactness}
%
%
The density of $\reps$ with respect to $\geps$, that is $\ueps = \reps/\geps$, is
smooth and satisfies the equation
\begin{align}
\label{uepseq}
\teps \geps \partial_t\ueps = \partial_x\bat{ \geps \partial_x\ueps }\,,
\end{align}
and hence we readily justify that
\begin{align}
\label{apriori1}
\int_{\Rset} \geps|\ueps|^2\,dx \,+\,\frac{2}{\teps} \int_0^t
\!\!\!\int_{\Rset} \geps|\partial_x \ueps|^2\,dx\,ds &=
\int_{\Rset} \geps|\ueps^0|^2\,dx ,
\\
\label{apriori2}
\frac{1}{2\teps}
\int_{\Rset} \geps|\partial_x \ueps|^2\,dx
\,+\, \int_0^t \!\!\!\int_{\Rset} \frac{|\partial_t \reps|^2}{\geps}\,dx\,ds
&= \frac{1}{2\teps}\int_{\Rset} \geps|\partial_x \ueps^0|^2\,dx.
\end{align}
Moreover, due to the assumption from Theorem \ref{T.1} the maximum principle
for \eqref{uepseq} implies that
\begin{align}%
\label{uepsbound}%
\inf\limits_{t\geq{0},\,x \in \R,\,\eps\geq0} \ueps\pair{t}{x}>0\,.%
\end{align}%
The a priori estimates \eqref{apriori1} and \eqref{apriori2} are direct consequences of the
$H^{-1}$-gradient and the $L^2$-gradient flow structures of \eqref{uepseq}, see \cite{PSV10}
for details. The Wasserstein structure, however, implies the a priori estimate
\begin{align}
\notag
\int_{\Rset}\reps\at{{\ln\reps}+\eps^{-2}H}\,dx
&+ \frac{1}{\teps} \int_0^t \!\!\!\int_{\Rset} \reps
\bat{\partial_x\at{\ln \reps + \eps^{-2}H}}^2\,dxds
= \int_{\Rset} \reps^0 \at{\ln \reps^0 + \eps^{-2}H}\,dx\,
\end{align}
and conserves the mass via $\int_{\Rset} \reps\,dx=1$. As mentioned before, our analysis does
not make use of this estimate but employs \eqref{apriori1} and \eqref{apriori2}.
%
%
\subsection{Rigorous formulation of the Rayleigh principle}\label{S.compactness.1}
We now derive a rigorous setting for the time integrated Rayleigh principle that corresponds
to the Wasserstein gradient structure of the Fokker-Planck equation. To this end we suppose
that $0<T<\infty$ is fixed and consider the weighted Lebesgue space
\begin{align*}
{L^2_\eps}:=\Big\{{f_\eps}:[0,T]\times\Rset\to\Rset\;:\;
\int_0^T\int_\Rset\frac{{f_\eps}^2}{\geps}\,dxdt<0\Big\},
\end{align*}
which is a Hilbert space for each $\eps>0$. We also define the linear space
\begin{align*}
\calT_{\eps}:=
\big\{%
\text{$\veps:L^2_{\eps}$ such that
$\veps=-\partial_x\at{\rho_\eps{j_\eps}}$ for some $j_\eps$ with
$\rho_\eps{j_\eps}\in{L^2_{\eps}}$}%
\big\},
\end{align*}
and show that the Rayleigh principle is a well-posed minimization problem on $\calT_\eps$.
\bigpar
\begin{lem}
\label{Lem:Tools}
For each ${f_\eps}\in{L^2_\eps}$ with $\partial_x{f_\eps}\in{L^2_\eps}$ we
have
\begin{align}
\label{Lem:Tools.Eqn1}
\int_0^T\int_\Rset\ln\ueps\partial_x{f_\eps} \, dxdt=-
\int_0^T\int_\Rset\frac{{f_\eps}\partial_x\ueps}{\ueps}\,dxdt,
\end{align}
and ${f_\eps}\pair{\cdot}{x}\to0$ as $x\to\pm\infty$ strongly in $L^2\bat{[0,T]}$. In particular,
we have
\begin{align*}
{f_\eps}\pair{t}{y}=\int_\infty^y\partial_x{f_\eps}\pair{t}{x}\,dx,\qquad
\int_\Rset\partial_x{f_\eps}\pair{t}{x}\,dx=0
\end{align*}
for almost all $\pair{t}{y}\in[0,T]\times\Rset$.
\end{lem}
\begin{proof}
Clearly, \eqref{Lem:Tools.Eqn1} holds for all smooth ${f_\eps}$ with compact support in
$[0,T]\times\Rset$. By approximation in ${L^2_\eps}$ -- and since we have
$\geps{\ln\ueps},\geps\ueps^{-1}\partial_x\ueps\in{L^2_\eps}$ due to \eqref{apriori1},
\eqref{apriori2}, \eqref{uepsbound}, and $\abs{\ln\ueps}\leq{C}\sqrt{\ueps}$ -- we then
conclude that \eqref{Lem:Tools.Eqn1} holds for all ${f_\eps}\in{L^2_\eps}$. Moreover, using
H\"older's inequality and $\int_{x_1}^{x_2}\geps\,dx\leq1$ we find
\begin{align*}
\int_0^T\babs{{f_\eps}\pair{t}{x_2}-{f_\eps}\pair{t}{x_1}}^2\,dt\leq
\int_0^T\abs{\int_{x_1}^{x_2}{\partial_x{f_\eps}\pair{t}{x}}\,dx}^2\,dt
\leq
\int_0^T\int_{x_1}^{x_2}\frac{\abs{\partial_x{f_\eps}}^2}{\geps}\,dxdt
\end{align*}
for all $x_1<x_2$. Therefore, and by assumption on $\partial_x{f_\eps}$, we know that
${f_\eps}\pair{\cdot}{x}$ converges strongly as $x\to\pm\infty$ to some limit functions in
$L^2\bat{[0,T]}$. From
\begin{align*}
\int_0^T\int_\Rset{f_\eps}^2\,dxdt\leq{\norm{\geps}_\infty}\int_0^T\int_\Rset\frac{{f_\eps}^2}{
\geps}\,dx\,dt<\infty
\end{align*}
we then infer that these limit functions vanish.
\end{proof}
\bigpar
\begin{cor}
\label{Cor:Tools}
We have $\partial_t\reps\in\calT_\eps$. Moreover, each
$v_\eps\in\calT_\eps$ satisfies
$\int_\Rset{v_\eps}\,dx=0$ for almost all $t$, as well as
\begin{align}
\label{Lem:Tools.Eqn2}
\abs{\int_0^T DE^{\eps}(\reps)\veps\,dt}^2\leq
\at{\int_0^T g^\eps_{\reps}\pair{\partial_t\reps}{\partial_t\reps}\,dt}
\at{\int_0^T g^\eps_{\reps}\pair{v_\eps}{v_\eps}\,dt}<\infty\,
\end{align}
with $g^\eps$ and $DE^\eps$ as in \eqref{eq12} and \eqref{eq14}.
\end{cor}
\begin{proof}
The first claim follows from $\partial_t\reps=\partial_x\at{\teps^{-1}\geps\partial_x\ueps}$
and since the a priori estimates \eqref{apriori1} and \eqref{apriori2} guarantee that
$\partial_t\reps\in{L^2_\eps}$ as well as $\teps^{-1}\geps\partial_x\ueps\in{L^2_\eps}$. Now
let $\veps=-\partial_xf_\eps\in\calT_\eps$ be arbitrary. By definition and \eqref{uepsbound},
we have
\begin{align*}
\int_0^T g^\eps_{\reps}\pair{v_\eps}{v_\eps}\,dt=
\int_0^T\int_\Rset\frac{\teps\,f_\eps^{\,2}}{\reps}\,dxdt\leq
{C}\int_0^T\int_\Rset\frac{\teps\,f_\eps^{\,2}}{\geps}\,dxdt<\infty\,.
\end{align*}
Moreover, Lemma \ref{Lem:Tools} provides $\int_\Rset{v_\eps}\,dx=0$ as well as
\begin{align*}
-\int_0^T DE^{\eps}(\reps)\veps\,dt
&=%
\int_0^T\int_{\Rset}f_\eps\partial_x\ln\ueps\,dxdt
=%
\int_0^T\int_\Rset
\frac{f_\eps}{\sqrt\reps}\frac{\geps\partial_x\ueps}{\sqrt\reps}\,dxdt
\\&\leq
\teps^2\at{\int_0^T\int_\Rset\frac{f_\eps^{\,2}}{\reps}\,dxdt}^{1/2}
\at{\int_0^T\int_\Rset\frac{\at{\teps^{-1}\geps\partial_x\ueps}^2}{\reps}\,dxdt}^{1/2},
\end{align*}
which is \eqref{Lem:Tools.Eqn2}.
\end{proof}
\bigpar
From Corollary \ref{Cor:Tools} we finally conclude that
\begin{align}
\label{Eqn:Integrated.Functional}
\int_0^{T} \tfrac 1 2 g^\eps_{\rho_\eps}\pair{\veps}{\veps} + DE^\eps(\reps)\veps\,dt
\end{align}
is well defined for $\veps\in\calT_\eps$, and that $\partial_t\reps$ is the unique minimizer.
%
%
\subsection{Compactness result for $\reps$}
%
We now exploit the a priori estimates \eqref{apriori1} and \eqref{apriori2} to derive
suitable compactness results for $\reps$, which then allow to extract convergent
subsequences. To this end we choose $0<\alpha<1$ independent of $\eps$ and 
define the intervals
\begin{align*}
J^\pm_\eps:=(\pm1-\eps^\alpha,\,\pm1+\eps^\alpha)\,,\qquad
J^0_\eps:=(-\eps^\alpha,\,+\eps^\alpha)\,,\qquad
\bar{J}_\eps:=\Rset\setminus\at{J^-_\eps{\cup}J^+_\eps}\,,
\end{align*}
as well as
\begin{align*}
I^-_\eps:=(-2+\eps^\alpha,\,-\eps^\alpha)\,,\qquad
I^+_\eps:=(\eps^\alpha,\,2-\eps^\alpha)\,,\qquad
I_{\eps}:=I^-_\eps\cup{I^+_\eps}\,.
\end{align*}
We proceed with summarizing some properties of $\geps$.
\bigpar
\begin{lem}%
\label{L.0}%
We have
\begin{align}%
\label{Lem:L.0.Eqn1}%
\int_{J^\pm_\eps}\geps\,dx
\quad\xrightarrow {\eps\to0}\quad\frac{1}{2}\,,
\qquad\qquad
\int_{\bar{J}_\eps}\geps\,dx
\quad\xrightarrow {\eps\to0}\quad0\,,
\end{align}%
as well as
\begin{align}%
\label{Lem:L.0.Eqn2}%
\sup\limits_{x\in\bar{J}_\eps}\geps\at{x}
\quad\xrightarrow {\eps\to0}\quad0\,,
\qquad\qquad
\tau_\eps^{-1}\inf_{x\in{I_\eps}}\geps\at{x}
\quad\xrightarrow{\eps\to0}\quad\infty\,,
\end{align}%
and%
\begin{align}%
\label{Lem:L.0.Eqn3}
\int_{J^0_\eps}
\frac{\teps}{\geps}\,dx\quad\xrightarrow{\eps\to0}\quad\frac{4}{k}\,,
\qquad\qquad
\int_{I_\eps}
\frac{\teps}{\geps}\,dx\quad\xrightarrow{\eps\to0}\quad0\,.
\end{align}
\end{lem}
\begin{figure}[t]%
\centering{%
\includegraphics[width=0.9\textwidth]{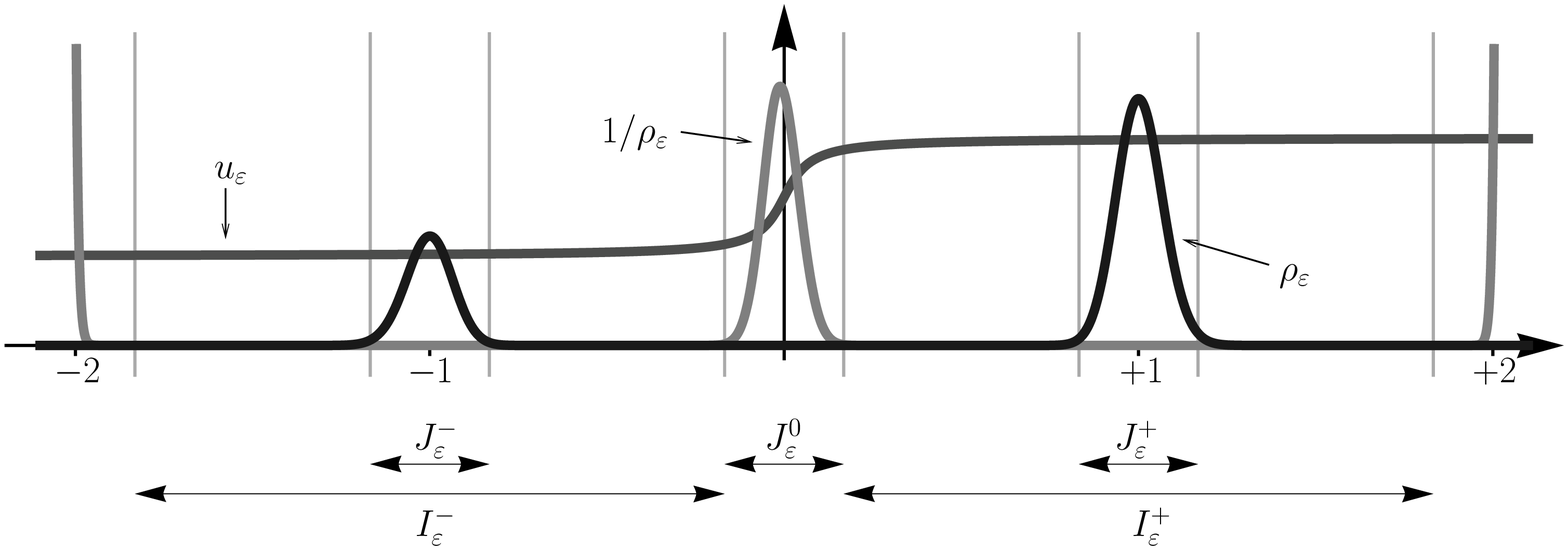}%
\vspace{-0.02\textheight}
}%
\caption{Schematic representation of $\reps$, $\ueps$, and $1/\reps$.}%
\label{FIG}%
\vspace{-0.03\textheight}
\end{figure}%
\begin{proof}
From the definition \eqref{gammadef} and $\int_\Rset\geps\,dx=1$ we readily derive
\eqref{Lem:L.0.Eqn1}. Thanks to our assumptions on $H$ we also have
\begin{align*}
H\at{x}\geq\min\bat{H\at{1+\eps^\alpha},\,H\at{1-\eps^\alpha}}\geq\tfrac{1}{2}
\at{1+o(1)}H^{\prime\prime}\at{1}\eps^{2\alpha}\qquad\forall\;
x\in\bar{J}_\eps\,,
\end{align*}
and this implies \eqref{Lem:L.0.Eqn2}$_1$.
Similarly, we find
\begin{align*}
H\at{x}\leq\max\bat{H\at{\eps^\alpha},\,H\at{2-\eps^\alpha}}\leq1+ o(1)-
\min\Bat{\tfrac{1}{2}H^{\prime\prime}\at{0}\eps^{2\alpha},\,H^\prime\at{2}\eps^\alpha}
\qquad\forall\;
x\in{I_\eps}
\end{align*}
and hence \eqref{Lem:L.0.Eqn2}$_2$ and \eqref{Lem:L.0.Eqn3}$_2$. Finally, a direct
computation gives
\begin{align*}
\int_{J^0_\eps} \frac{\teps}{\geps}\,dx
&=  \frac{Z_{\eps}}{\eps^2}\int_{-\eps^\alpha}^{\eps^\alpha}\exp\at{-
\tfrac{1}{2}H^{\prime\prime}\at{0}\bat{1+o(1)}x^2/\eps^{2} } \,dx
= \frac{ 4 \pi\,\bat{1+o(1)}}{\sqrt{ H^{''}(1) |H^{''}(0)|}}
\end{align*}
which is \eqref{Lem:L.0.Eqn3}$_1$ thanks to \eqref{kramers}.
\end{proof}
\bigpar
Our compactness result for $\reps$ is illustrated in Figure \ref{FIG} and reads as follows.
\bigpar%
\begin{lem}
\label{Lem:Compactness}
There exists a subsequence $\eps\to0$ along with a function $u\in{H^1\at{[0,T]}}$ with weak
derivative $\dot{u}$ such that
\begin{enumerate}
\item

for each $t\in[0,T]$, $\reps\pair{t}{\cdot}$ is a probability measure on $\Rset$ that
converges weakly$\star$ to
$\tfrac{1}{2}\bat{u\at{t}\delta_{+1}+\at{1-u(t)}\delta_{-1}}$ such that
\begin{align}
\label{Lem:Compactness.Eqn21}
\int_{J^\pm_\eps}\reps\,dx\quad\xrightarrow{\eps\to0}\quad
\tfrac{1}{2}\pm\tfrac{1}{2}\at{u-1}
\qquad\text{and}\qquad
\int_{\bar{J}_\eps}\reps\,dx
\quad\xrightarrow{\eps\to0}\quad0\,
\end{align}
uniformly in $t\in[0,T]$,
\item
$\partial_t\reps$ is a measure on $[0,T]\times\Rset$ that converges weakly$\star$ to
$\tfrac{1}{2}\bat{\dot{u}\delta_{+1}-\dot{u}\delta_{-1}}$ such that
\begin{align}
\label{Lem:Compactness.Eqn10}
\int_{J^\pm_\eps}\partial_t\reps\,dx\quad\xrightharpoonup{\eps\to0}\quad
\pm\tfrac{1}{2}\dot{u}
\qquad\text{and}\qquad
\int_{\bar{J}_\eps}\abs{\partial_t\reps}\,dx\quad\xrightarrow{\eps\to0}\quad0
\end{align}
weakly and strongly in $L^2\bat{[0,T]}$, respectively.
\end{enumerate}
Moreover, we have
\begin{math}
c\;\leq \;u\at{t},\,2-u\at{t}\;\leq\;2-c
\end{math} %
for some $0<c<2$ and all $t\in[0,T]$.
\end{lem}
\begin{proof}
The estimates \eqref{apriori1} and \eqref{apriori2} combined with \eqref{Lem:L.0.Eqn2} imply
\begin{align*}
\sup_{t\in[0,T]}\at{\int_{\bar{J}_\eps}\reps\,dx}^2
&\leq%
\sup_{t\in[0,T]}\at{\int_{\Rset}\geps\ueps^2\,dx}
\at{\int_{\bar{J}_\eps}\geps\,dx}
\quad\xrightarrow{\eps\to0}\quad0\,,
\\%
\int_0^T\at{\int_{\bar{J}_\eps}\abs{\partial_t\reps}dx}^2\,dt
&\leq%
\int_0^T\at{\int_{\bar{J}_\eps}
\geps\abs{\partial_t\ueps}^2\,dx }
\at{\int_{\bar{J}_\eps}\geps\,dx}\,dt
\quad\xrightarrow{\eps\to0}\quad0\,,
\end{align*}
and hence we find \eqref{Lem:Compactness.Eqn21}$_2$ and \eqref{Lem:Compactness.Eqn10}$_2$. We
now consider the functions $u^\pm_\eps:[0,T]\to\Rset$ defined by
\begin{align*}
\tfrac{1}{2}u^\pm_\eps\at{t}
:=
\int_{J^\pm_\eps}\reps\pair{t}{x}\,dx=
\int_{J^\pm_\eps}\ueps\pair{t}{x}\geps\at{x}\,dx\,,
\end{align*}
and using \eqref{apriori1}, \eqref{apriori2}, and \eqref{Lem:L.0.Eqn2} again we obtain
\begin{align*}
\tfrac{1}{4}\norm{u^\pm_\eps}^2_{L^2\at{[0,T]}}&\leq\int_0^T\at{\int_{J^\pm_\eps}\geps\ueps^2\,
dx}
\at{\int_{J^\pm_\eps}\geps\,dx}\,dt\leq{C}\,,
\\%
\tfrac{1}{4}\norm{\dot{u}^\pm_\eps}^2_{L^2\at{[0,T]}}&\leq\int_0^T\at{\int_{J^\pm_\eps}
\geps\abs{\partial_t\ueps}^2\,dx }
\at{\int_{J^\pm_\eps}\geps\,dx}\,dt\leq{C}\,.
\end{align*}
By weak compactness we can therefore extract a subsequence such that
\begin{align*}
u^\pm_\eps\;\xrightharpoonup{\eps\to0}\;u^\pm\quad\text{weakly in
$H^1\bat{[0,T]}\subset\subset{C\bat{[0,T]}}$},\quad
\end{align*}
for some limit functions $u^\pm$. Setting $u:=u_+$, the remaining assertions hold either by
construction, or thanks to \eqref{uepsbound} and $\int\reps\,dx=1$ .
\end{proof}
\bigpar
From now on we suppose that a subsequence as in Lemma \ref{Lem:Compactness} is chosen and
prove the assertions of Theorem \ref{T.1} for this subsequence. Since the function $u$ is
uniquely determined by the limit problem \eqref{limiteq}, we then conclude afterwards
that Theorem \ref{T.1} holds for all subsequences.
%
%
\subsection{Leading order description of $\ueps$}
%
In order to pass to the limit $\eps\to\infty$ in the minimization problem corresponding to
\eqref{Eqn:Integrated.Functional} we show that the relative density $u_\eps$ is close to a
step function but exhibits a narrow and smooth transition layer at $x=0$ whose shape is
determined by $\geps$, see Figure \ref{FIG}. Specifically, we prove that $u_\eps$ can be
approximated by
\begin{align}
\label{Approx.Eqn1}
\tilde{u}_\eps\pair{t}{x} := 1+\bat{u\at{t}-1}\,\eta_\eps\at{x}
\,,\qquad
\eta_\eps\at{x} := \displaystyle2\int_0^x
\frac{dy}{\geps\at{y}}\;\Big/\;\int_{-1}^1\frac{dy}{\geps\at{y}}\,.
\end{align}
Notice that $\eta_\eps$ is an approximation of the sign function on the interval
$\abs{x}<2$ because Lemma \ref{L.0} provides that 
\begin{align}
\label{Eqn:Props.Eta}
\sup_{x\in{I_\eps}}\abs{\eta_\eps\at{x}-\mathrm{sign}\at{x}}\xrightarrow{\eps\to0}0.
\end{align} %

\begin{lem}
\label{L.1a}
We have
\begin{align*}
\sup_{t\in[0,T],x\in{I^\pm_\eps}}
\babs{\ueps\pair{t}{x}-\tilde{u}_\eps\pair{t}{x}}
\;\xrightarrow{\eps\to0}\;0%
\,,\qquad
\int_0^{T} \sup_{ x \in J^0_{\eps}}
\Big | \ueps(t,x) - \tilde{u}_\eps\pair{t}{x} \Big|^2 \,dt
\;\xrightarrow{\eps\to0}\;0\,.
\end{align*}
\end{lem}
\begin{proof}
From \eqref{apriori2} and \eqref{Lem:L.0.Eqn2}$_2$ we infer that
\begin{align*}
\sup_{t\in[0,T]}\int_{I_\eps}\abs{\partial_x\ueps}^2\,dx
\leq
\sup_{t\in[0,T]}\at{\int_{I_\eps}\frac{\geps\abs{\partial_x\ueps}^2}{\teps}\,dx}
\at{\int_{I_\eps}\frac{\teps}{\geps}\,dx}
\xrightarrow{\eps\to0}0\,,
\end{align*}
and hence \begin{math}
\sup_{t\in[0,T]}\sup_{x,y\in{I^\pm_\eps}}\babs{\ueps\pair{t}{x}-
\ueps\pair{t}{y}}\xrightarrow{\eps\to0}0.
\end{math} %
With \eqref{Lem:L.0.Eqn1} we therefore conclude that
\begin{align*}
\sup\limits_{t\in[0,T]}\sup\limits_{x\in{I^\pm_\eps}}\abs{\ueps\pair{t}{x}-2\int_{J^\pm_\eps}
\reps\pair{t}{y}\,dy}
\quad\xrightarrow{\eps\to0}\quad0\,,
\end{align*}
so Lemma \ref{Lem:Compactness} gives
\begin{align*}
\sup\limits_{t\in[0,T]}\sup\limits_{x\in{I^\pm_\eps}}\babs{1\pm\at{u-1}-\ueps\pair{t}{x}}
\quad\xrightarrow{\eps\to0}\quad0\,.
\end{align*}
The first assertion now follows from \eqref{Eqn:Props.Eta}. Towards the second claim we
integrate \eqref{uepseq} twice with respect to $x$. This gives
\begin{align}
\label{L.1a.Eqn1}
\ueps\pair{t}{x} &= \hat{u}_\eps\pair{t}{x} +
\int_0^x \frac{\teps}{\geps\at{y}} \int_0^y
\partial_t \reps\pair{t}{z} \,dz \,dy
\end{align}
with
\begin{math}
\hat{u}_\eps\pair{t}{x}:=
C_{1,\eps}\at{t} \eta_\eps\at{x} +
C_{2,\eps}\at{t},
\end{math} %
where the two constants of integration can be computed by
\begin{align}
\label{L.1a.Eqn11}
C_{1,\eps}\at{t}=
\frac{\hat{u}_\eps\pair{t}{+1/2}-\hat{u}_\eps\pair{t}{-1/2}}{2\eta_\eps\at{1/2}},\qquad
C_{2,\eps}\at{t}=
\frac{\hat{u}_\eps\pair{t}{+1/2}+\hat{u}_\eps\pair{t}{-1/2}}{2}.
\end{align}
For
$\abs{y}\leq\abs{x}\leq1-\eps^\alpha$ the integral term in
\eqref{L.1a.Eqn1} can be estimated by
\begin{align*}
\abs{\int_0^{y} \partial_t \reps\pair{t}{z}\,dz}&\leq \Big( \int_0^{1-\eps^\alpha}
\frac{|\partial_t
\reps\pair{t}{z}|^2}{\geps\at{z}}\,dz
\Big)^{1/2}
\Big( \int_0^{1-\eps^\alpha} \geps\at{z} \,dz \Big)^{1/2}
\\& \leq C
\at{\sup\limits_{\abs{x}\leq{1}-\eps^\alpha}\geps\at{x}}^{1/2}
\at{\int_{\Rset} \frac{|\partial_t \reps\pair{t}{z}|^2}{\geps\at{z}} dz}^{1/2}\,,
\end{align*}
so \eqref{apriori2} and Lemma \ref{L.0} imply
\begin{align*}
\int_0^{T}\sup_{\abs{x}\leq{1}-\eps^\alpha}
\Big | \ueps(t,x) -\hat{u}_\eps\pair{t}{x}\Big|^2 \,dt
\quad\xrightarrow{\eps\to0 }\quad0\,.
\end{align*}
In particular, we have
$\hat{u}_\eps\pair{\cdot}{\pm1/2}
\xrightarrow{\eps\to0}1\pm\at{u-1}$
in $L^2\bat{[0,T]}$. From this, \eqref{L.1a.Eqn11}, and the definition of both $\hat{u}_\eps$
and $\tilde{u}_\eps$ we finally conclude that
\begin{align*}
\int_0^{T} \sup_{\abs{x}\leq{2}-\eps^\alpha}
\Big | \hat{u}_\eps\pair{t}{x} -\tilde{u}_\eps\pair{t}{x}\Big|^2 \,dt
\quad\xrightarrow{\eps\to0}\quad0\,,
\end{align*}
and the proof is complete.
\end{proof}
\bigpar
It follows from Lemma \ref{L.0} that for $\eps\to0$ the functions $\tau_\eps/\geps$ generate
a delta distribution in $x=0$ with height $4/k$. The following result combined with Lemma
\ref{L.1a} shows that $\tau_\eps/\reps$ has a similar property. For the proof we 
recall that the function $u\mapsto\frac{1}{u-1}\ln\frac{u}{2-u}$ is continous in
$u=1$ and uniformly positive for $u\in(0,2)$.
\bigpar
\begin{lem}
\label{L.2}
We have
\begin{math}
\displaystyle%
\int_{J^0_{\eps}}\frac{\teps}{\geps\tilde{u}_\eps}\,dx \;\xrightarrow{\eps\to0}\,
\frac{2}{k\at{u-1}}\ln{\frac{u}{2-u}}
\end{math}
uniformly in $t\in[0,T]$.
\end{lem}
\begin{proof}
By definition we have
\begin{align*}
\frac{\teps}{\geps\at{x}\tilde{u}_\eps\pair{t}{x}}=
\frac{C_\eps}{2\at{u\at{t}-1}}\,\partial_x\bat{\ln\tilde{u}_\eps\pair{t}{x}}
\,,\qquad
C_\eps^{-1}:=\teps\int_{-1}^{-1}\frac{1}{\geps}\,dx\,.
\end{align*}
This implies
\begin{align*}
\int_{J^0_{\eps}}\frac{\teps}{\geps\at{x}\tilde{u}_\eps\pair{t}{x}}\,dx=
\frac{C_\eps}{2\at{u\at{t}-1}}
\ln\frac{\tilde{u}_\eps\pair{t}{+\eps^\alpha}}{\tilde{u}_\eps\pair{t}{-\eps^\alpha}}
\,,
\end{align*}
so the desired result follows from \eqref{Approx.Eqn1} and \eqref{Eqn:Props.Eta}.
\end{proof}
\bigpar
%
%
\section{Passing to the limit in the Rayleigh principle}\label{S.proofs}
%
%
We are now able to prove our main result, which then implies Theorem \ref{T.1}.
\bigpar
\begin{thm}
\notag
The limit of $\partial_t\reps$, that is $\dot{u}$, satisfies
\begin{align*}
\int_0^{T} \Big ( \tfrac 1 2 g_{u}\big(\dot u, \dot u\big)
+ DE(u) \dot u \Big)\,dt
\leq \int_0^{T}  \Big ( \tfrac 1 2 g_{u}\big(v, v\big)
+ DE(u) v \Big)\,dt
\end{align*}
for all $v\in L^2\bat{[0,T]}$.
\end{thm}
\bigpar
This theorem is implied by the following three Lemmas, which ensure that the limit $\dot u$
of the minimizers $\partial_t \reps$ of the Rayleigh principle associated to the
$\eps$-problems is indeed a minimizer of the Rayleigh principle associated to the limit
problem. Notice that the assertions of Lemmas \ref{L.3}-\ref{L.5} are closely related, but
not equivalent to the $\Gamma$-convergence of the Rayleigh principles. In fact, in Lemma
\ref{L.3} we prove lower semi-continuity of the metric tensor only for the minimizers
$\partial_t \reps$. This is, however, sufficient to conclude that $\dot u$ is again a
minimizer of the limit problem.
\bigpar%
\begin{lem}(lim-inf estimate for metric tensor)
\label{L.3}
We have
\begin{align*}
\liminf_{\eps \to 0} \int_0^{T}
g^{\eps}_{\reps} \big( \partial_t \reps, \partial_t \reps\big)\,dt
\geq \int_0^{T}
g_{u}\big(\dot u, \dot u \big)\,dt\,.
\end{align*}
\end{lem}
\begin{proof}
Recall that we have
\begin{align*}
\int_0^{T} g_{\reps}^{\eps} \big( \partial_t \reps, \partial_t \reps
\big)\,dt
=
\int_0^{T} \int_{\Rset} \frac{\teps f_{\eps}^2}{\reps}\,dx\,dt
\geq
\int_0^{T} \int_{J^0_\eps}\frac{\teps f_{\eps}^2}{\reps}\,dx\,dt\,,
\end{align*}
where $f_\eps:=\teps^{-1}\geps\partial_x\ueps$ satisfies the a priori estimate
\begin{align*}
\int_{J^0_\eps}\frac{\teps\,f_\eps^2}{\reps}\,dx=
\int_{J^0_\eps}\frac{\geps\,\abs{\partial_x\ueps}^2}{\teps\ueps}\,dx
\leq%
{C}\quad\text{for all
$t\in[0,T]$}\,.
\end{align*}
Consequently, Lemma \ref{L.1a} and \eqref{uepsbound} ensure that
\begin{align*}
\int_{J^0_\eps}\frac{\teps\,f_\eps^2}{\reps}-
\int_{J^0_\eps}\frac{\teps\,f_\eps^2}{\geps\tilde{u}_\eps}
\quad\xrightarrow{\eps\to0}\quad
0
\quad\text{strongly in $L^2\bat{[0,T]}$}\,,
\end{align*}
and therefore it is sufficient to show that
\begin{align*}
\liminf\limits_{\eps\to0}\int_0^{T} \int_{J^0_\eps}
\frac{\teps f_{\eps}^2}{\geps\tilde{u}_\eps}\,dx\,dt
\geq
\int_0^{T} g_{u}\big(\dot u, \dot u \big)\,dt\,.
\end{align*}
Setting $\tilde{f}_\eps\at{t}:={f}_\eps\pair{t}{0}$ and using
$\partial_x{f_\eps}=\partial_t\reps$ we find
\begin{align*}
\sup\limits_{x\in{J^0_\eps}}%
\babs{\tilde{f}_\eps\at{t}-f_\eps\pair{t}{x}}\leq
\at{\int_{x\in{J^0_\eps}}\frac{\abs{\partial_t\reps}^2}{\geps}
\,dx}^{1/2}\at{\int_{J^0_\eps}{\geps}\,dx}^{1/2},
\end{align*}
and hence
\begin{align*}
\abs{\int_{J^0_\eps}\frac{\teps{f_\eps}}{\geps\tilde{u}_\eps}-
\int_{J^0_\eps}\frac{\teps{\tilde{f}_\eps}}{\geps\tilde{u}_\eps}}^2
\leq
\at{\int_{x\in{J^0_\eps}}\frac{\abs{\partial_t\reps}^2}{\geps}
\,dx}\at{\int_{J^0_\eps}{\geps}\,dx}\at{\int_{J^0_\eps}\frac{\teps}{\geps\tilde{u}_{\eps}}\,dx
}^2.
\end{align*}
Integration with respect to $t$, and employing \eqref{apriori2} as well as Lemma \ref{L.0},
then gives
\begin{align}
\label{L.3.Eqn1}
\int_{J^0_\eps}\frac{\teps{f_\eps}}{\geps\tilde{u}_\eps}\;-
{\tilde{f}_\eps}\int_{J^0_\eps}\frac{\teps}{\geps\tilde{u}_\eps}
\quad\xrightarrow{\eps\to0}\quad0
\quad\text{strongly in $L^2\bat{[0,T]}$}\,.
\end{align}
Applying Jensen's inequality, the convergence \eqref{L.3.Eqn1}, and Lemma \ref{L.2} we 
estimate
\begin{align*}
\liminf\limits_{\eps\to0}
\int_0^T\int_{J^0_\eps}
\frac{\teps f_{\eps}^2}{\geps\tilde{u}_\eps}\,dx\,dt
&\geq
\liminf\limits_{\eps\to0}
\int_0^T
\at{\int_{J^0_\eps}
\frac{\teps f_{\eps}}{\geps\tilde{u}_\eps}\,dx}^2
\at{\int_{J^0_\eps}
\frac{\teps}{\geps\tilde{u}_\eps}\,dx}^{-1}\,dt
\\&\geq
\liminf\limits_{\eps\to0}
\int_0^T {\tilde{f}_{\eps}}^{\,2}
\at{\int_{J^0_\eps}
\frac{\teps}{\geps\tilde{u}_\eps}\,dx}\,dt.
\end{align*}
The desired result now follows from Lemma \ref{L.2} and since
Lemma \ref{Lem:Tools}, combined with Corollary \ref{Cor:Tools} and
\eqref{Lem:Compactness.Eqn10}, implies that
$\tilde{f}_\eps=\int_\infty^0\partial_t\reps\,dx\xrightharpoonup{\eps\to0}\tfrac{1}{2}\dot{u}$
weakly in
$L^2\bat{[0,T]}$.
\end{proof}
\bigpar
\begin{lem}(lim-inf estimate for energy)
\label{L.4}
We have
\begin{align*}
\int_0^{T} DE^{\eps}(\reps) \partial_t \reps \,dt
\quad\xrightarrow{\eps\to0}\quad \int_0^{T} DE(u) \dot u \,dt\,.
\end{align*}
\end{lem}%
\begin{proof}
By \eqref{uepsbound}, which implies that $\abs{\ln\ueps}\leq{C}\sqrt{\ueps}$, and using
\eqref{apriori1}, \eqref{apriori2}, and \eqref{Lem:L.0.Eqn1}, we find
\begin{align*}
\Big|\int_0^{T} \int_{\bar{J}_\eps} \partial_t \reps \ln \ueps
\,dx\,dt\Big|^2
&\leq  C\at{\int_0^{T} \int_{\bar{J}_\eps} \frac{|\partial_t \reps|^2}{\geps}
\,dx\,dt}
\at{\int_0^{T} \int_{\bar{J}_\eps}\geps{\ueps}\,dx\,dt}
\\&\leq C \sup_{0\leq{t}\leq{T}} \at{\int_{\bar{J}_\eps}  \geps|\ueps|^2
\,dx}^{1/2} \at{
\int_{\bar{J}_\eps} \geps \,dx }^{1/2}
\quad\xrightarrow{\eps\to 0}\quad0\,.
\end{align*}
Moreover, the convergence results from Lemma \ref{Lem:Compactness} and Lemma \ref{L.1a}
provide
\begin{align*}
\int_0^{T} \int_{J^\pm_\eps} \partial_t \reps \at{\ln \ueps} \,dx \,dt
\quad\xrightarrow{\eps\to0}\quad
\int_0^{T} \frac{\dot{u}}{2}\,\ln{\frac{u}{2-u}}\,dt\,,
\end{align*}
which is the desired result.
\end{proof}
\bigpar
\begin{lem}(Existence of recovery sequence)
\label{L.5}
For all $v \in L^2\bat{[0,T]}$ there exists a sequence $\veps\in\calT_{\eps}$
 such that
\begin{align*}
\limsup_{\eps \to 0} \int_0^{T} \tfrac 1 2 g^{\eps}_{\reps} \big( \veps,
\veps\big) + DE^{\eps}(\reps) \veps\,dt \leq \int_0^{T} \tfrac 1 2
g_{u}\big(v,v\big) +  DE(u) v\,dt\,.
\end{align*}
\end{lem}
\begin{proof}
For each $\eps>0$ we define $\veps\pair{t}{x}:=\tfrac1 2\,v(t)\big( \psi_{{\eps}}(x-1) -
\psi_{{\eps}}(x+1)\big)$, where $\psi_{{\eps}}$ is a standard Dirac sequence with support in
$(-{\eps^\alpha},{\eps^\alpha})$. We also introduce functions $f_\eps$ by $f_\eps\pair{t}{x}
:= - \int_{-\infty}^x \veps(t,y) \,dy$. By construction, $f_\eps$ vanishes for
$\abs{x}\geq{1+\eps^\alpha}$, is equal to $\frac12v(t)$ for $\abs{x}\leq1-\eps^\alpha$, and
satisfies
\begin{align*}
0\leq\int_0^{T} \int_{J^\pm_\eps} \frac{\teps f_\eps^{\,2}}{\reps}\,dx\,dt
\leq\int_0^{T} \tfrac{1}{4}v^2\int_{J^\pm_\eps} \frac{\teps }{\geps\ueps}\,dx\,dt
\quad\xrightarrow{\eps\to 0}\quad0
\end{align*}
thanks to \eqref{uepsbound} and \eqref{Lem:L.0.Eqn3}. Moreover, in view of
Lemma \ref{L.1a}, Lemma \ref{L.2} and \eqref{Lem:L.0.Eqn3} we also have
\begin{align*}
\int_0^{T} \int_{\bar{J}_\eps} \frac{\teps f_\eps^{\,2}}{\reps}\,dx\,dt
=\int_0^{T} \tfrac{1}{4}v^2\int_{\abs{x}\leq{1}-\eps^\alpha} \frac{\teps }{\geps\ueps}\,dx\,dt
\quad\xrightarrow{\eps\to 0}\quad
\int_0^T\frac{v^2}{2k\at{u-1}}\ln{\frac{u}{2-u}}\,dt,\,
\end{align*}
and hence
\begin{align*}
\int_0^{T} g^{\eps}_{\reps} \big( \veps,
\veps\big) \,dt =
\int_0^{T} \int_{\Rset} \frac{\teps f_\eps^{\,2}}{\reps}\,dx\,dt
\quad\xrightarrow{\eps\to 0}\quad
\int_0^{T}
g_{u}\big(v,v\big)  v\,dt\,.
\end{align*}
Finally,
\begin{math}
\int_0^{T} \int_{\Rset} \veps \ln \ueps \,dx \,dt
\xrightarrow{\eps\to0}\int_0^T\frac{v}{2}\ln \frac{u}{2-u} \,dt
\end{math}
follows as in the proof of Lemma \ref{L.4}.
\end{proof}%
\section*{Acknowledgements} The authors gratefully to Alexander Mielke for stimulating
discussions and to the referee for his valuable comments.
This work was supported by the EPSRC Science and Innovation award to the Oxford Centre for
Nonlinear PDE (EP/E035027/1).

\end{document}